\documentclass{IEEEtran4PSCC}
\ifCLASSINFOpdf
   \usepackage[pdftex]{graphicx}
\else
   \usepackage[dvips]{graphicx}
\fi
\usepackage[cmex10]{amsmath}
\usepackage{color,verbatim}
\usepackage{siunitx}
\usepackage[shortlabels]{enumitem}
\usepackage{bbm}
\usepackage{amsmath,amssymb}
\usepackage[ruled,vlined,linesnumbered]{algorithm2e}

\usepackage{amssymb}
\usepackage{amsmath}
\usepackage{amsthm}
\usepackage[english]{babel}
\usepackage{xspace}
\usepackage{cases}
\usepackage{color}
\usepackage{graphicx}
\usepackage[caption=false]{subfig}
\usepackage{epstopdf}
\usepackage{upgreek}
\usepackage{booktabs}
\usepackage{array}
\usepackage{multirow}
\usepackage{bbold}
\usepackage[hidelinks]{hyperref}
\usepackage{soul}
\usepackage{enumerate}
\usepackage{enumitem}
\usepackage{csquotes}
\usepackage{subfloat}

\newcommand{\vplain}{v}
\newcommand{\vre}{\vplain^{\text{re}}}

\newcommand{\vrebold}{\boldsymbol{\vplain}^{\mathrm{re}}}
\newcommand{\vim}{\vplain^{\text{im}}}%

\newcommand{\vimbold}{\boldsymbol{\vplain}^{\mathrm{im}}}

\newcommand{\pce}{\textsc{pce}\xspace}

\newcommand{\xib}{\boldsymbol{\xi}}
\newcommand{\Psib}{\mathbf{\Psi}}

\newcommand{\TV}{\operatorname{TV}}

\DeclareMathOperator{\E}{\mathbb{E}}

\newcommand{\Co}{C_{\rm{off}}}
\newcommand{\algo}{\textsc{Spice}}

\usepackage{tikz}
\def\checkmark{\tikz\fill[scale=0.4](0,.35) -- (.25,0) -- (1,.7) -- (.25,.15) -- cycle;}

\allowdisplaybreaks

\makeatletter
\let\old@ps@headings\ps@headings
\let\old@ps@IEEEtitlepagestyle\ps@IEEEtitlepagestyle
\def\psccfooter#1{%
    \def\ps@headings{%
        \old@ps@headings%
        \def\@oddfoot{\strut\hfill#1\hfill\strut}%
        \def\@evenfoot{\strut\hfill#1\hfill\strut}%
    }%
    \def\ps@IEEEtitlepagestyle{%
        \old@ps@IEEEtitlepagestyle%
        \def\@oddfoot{\strut\hfill#1\hfill\strut}%
        \def\@evenfoot{\strut\hfill#1\hfill\strut}%
    }%
    \ps@headings%
}
\makeatother

\begin{document}
%

\title{Efficient Polynomial Chaos Expansion for Uncertainty Quantification in Power Systems}

\author{
\IEEEauthorblockN{David M\'etivier, Marc Vuffray, Sidhant Misra}
\IEEEauthorblockA{Los Alamos National Laboratory \\
\{metivier, vuffray, sidhant\}@lanl.gov} }

\maketitle



\maketitle

\begin{abstract}
Growing uncertainty from renewable energy integration and distributed energy resources motivate the need for
advanced tools to quantify the effect of uncertainty and assess the risks it poses to secure system operation. Polynomial chaos expansion (PCE) has been recently proposed as a tool for uncertainty quantification in power systems. The method produces results that are highly accurate, but has proved to be computationally challenging to scale to large systems. We propose a modified algorithm based on PCE with significantly improved computational efficiency that retains the desired high level of accuracy of the standard PCE. Our method uses computational enhancements by exploiting the sparsity structure and algebraic properties of the power flow equations. We show the scalability of the method on the 1354 pegase test system, assess the quality of the uncertainty quantification in terms of accuracy and robustness, and demonstrate an example application to solving the chance constrained optimal power flow problem.

\end{abstract}


\thanksto{\noindent Submitted to the 21st Power Systems Computation Conference (PSCC 2020).}

\section{Introduction}
Traditional power systems operational planning and management is being challenged by the increased penetration of renewable energy and distributed energy resources. The variability of power consumption and generation inherent to these additions calls for new control and optimization tools capable of accurately handling the impact of uncertainty on a faithful, nonlinear description of the power grid. However, the  non-linearity introduces significant computational challenges in quantifying the effect of uncertainty on the system, spurring a long line of research on the topic. 

The most commonly used approach is based on the linear DC approximation approximation to the AC power flow equations (AC-PFE) \cite{Bienstock14}, \cite{roald2016corrective}. The algebraic simplicity facilitates fast computations at the cost of accuracy, which can be significant when the uncertainties are large. More accurate approaches \cite{Roald18} use a hybrid representation, where the full non-linear equations are used for the nominal power flows and the effect of uncertainty is linearized, and are appropriate for moderate uncertainty magnitudes. In contrast, methods based on Monte Carlo are accurate and capture the non-linear implicit nature of the AC power flow equations in a faithful way. But attaining a sufficient precision using Monte Carlo involves solving the system of equations for a large number of random uncertainty realizations which can result in unacceptably large computation times.

A recent line of work \cite{Muehlpfordt16b,Muehlpfordt18b,Muehlpfordt17a,Muehlpfordt2019TPWRS} proposes the use of polynomial chaos expansion (PCE) to handle the non-linear AC-PFE. Using PCE, all uncertain quantities in the system are expressed as polynomials of the uncertain variables. The coefficients of the polynomial are tailored to the uncertainty distribution by performing an orthogonal projection step. Uncertainty quantification with PCE reduces to solving an extended system of power flow like equations, which we call the PCE overloaded power flow, and the method lends itself to easy integration into uncertainty-aware economic dispatch problems such as the chance constrained optimal power flow. However in its current form, PCE lacks sufficient scalability that precludes its use for large power systems.

 In this work, we develop a PCE-based method coined \algo{} for Sparse Polynomial Iterative Chaos Expansion. The iterative procedure in \algo{} identifies and exploits the sparsity structure inherent to the topology of the grid and reflected in the PCE overloaded power flow equations. Approximations based on the algebraic properties of the power flow equations are used to further simplify the problem. By a careful employment of the simplifications, \algo{} is able to significantly reduce the computational complexity of standard PCE, while still retaining its accuracy. 
 
 We demonstrate the improvements in scalability with a  detailed numerical study on the 1354 bus pegase test system. We show that the polynomials produced by \algo{} can be used to perform highly accurate uncertainty quantification, and therefore Monte Carlo simulations can be carried out without the need to repeatedly solve the power flow equations. As an application, we use the iterative procedure developed in \cite{Roald18} to solve the chance constraint optimal power flow problem.

\section{Power Flow equations}\label{sec:PFE}
The power network is modeled using a graph with $N$ buses and $L$ transmission lines. We use $(p_i,q_i)$ to denote the net active and reactive power injection at bus $i$. The power flow physics is described by a system of quadratic equations known as the Kirchoff's laws and are given by 
\begin{subequations}    \label{eq:pfe}
\begin{align}   
    p_i &= \sum_{j \in N} G_{ij} ( \vre_{i} \vre_{j} + \vim_{i} \vim_{j} ) + B_{ij} ( \vim_{i} \vim_{j} - \vre_{i} \vim_{j} ),\\
q_i &= \sum_{j \in N} G_{ij} ( \vim_{i} \vim_{j} - \vre_{i} \vim_{j} ) - B_{ij} ( \vre_{i} \vre_{j} + \vim_{i} \vim_{j} ),
\end{align}
\end{subequations}

where $G_{ij} + j B_{ij}$ denotes the $(ij)^{th}$ entry of the complex bus impedance matrix, and $\vre_{i}$ and $\vim_{i}$ denote the real and imaginary part of the complex voltage phasor at bus $i$. The equations in \eqref{eq:pfe} are known as the AC power flow equations (AC-PFE) in rectangular coordinates. In a more abstract form, the power flow equations are a system of $2N$ quadratic equations that map the voltage phasors to the bus injections. We denote these by 
\begin{subequations}    \label{eq:pfe_compact}
\begin{align}
    p_i = p_i(\vrebold,\vimbold), \quad q_i = q_i(\vrebold,\vimbold),
\end{align}
\end{subequations}
where $p_i(), q_i()$ are the quadratic functions described in \eqref{eq:pfe}.
The non-linear nature of the AC-PFE are a significant mathematical challenge in many problems important to power systems planning and operations. These include (i) the \emph{AC optimal power flow problem (AC-OPF)} used for economic generation dispatch where they appear as non-convex constraints, and (ii) \emph{uncertainty quantification (UQ)} used to analyze the effect of uncertainty in the net power injection at buses $p_i,q_i$ caused by load and renewable generation.

The latter is particularly challenging, especially since the uncertainty quantification methods are often required to be incorporated within an optimization framework such as stochastic and robust variants of the AC-OPF. In this paper, we aim at developing an uncertainty quantification method that is both scalable and accurate. We adopt the recently proposed approach based on the so-called polynomial chaos expansion. While PCE has been shown to enable compact and accurate UQ, it suffers from the curse of dimensionality described in the next section.


\section{Polynomal Chaos Expansion}
In this section, we provide a brief overview of the polynomial chaos expansion approach. For a detailed exposition, the reader is refered to \cite{Muehlpfordt16b}. 

When the power system is subject to uncertain power injections, i.e., when the quantities $p_i,q_i$ in \eqref{eq:pfe}
are random variables, this randomness propagates through the system of equations resulting in every other variables (voltages, line power flows, etc.) behaving as random variables. The state of the system is therefore a \emph{function} of the uncertainty realization making it inherently difficult to obtain a compact representation of the system behavior. Polynomial chaos deals with this problem by using using a polynomial representation for each of these functions. Further, instead of using the standard monomial basis, PCE uses a special set of basis functions for the polynomial expansion that are \emph{orthogonal} with respect to the uncertainty distribution. 

\subsection{Uncertainty Quantification Using PCE}    \label{subsec:uq_pce}
Let $\xib=(\xi_1,\cdots,\xi_n)$ denote the vector of random variables, where $n$ is the dimension of the uncertainty. These variables may be used to directly represent the random variables corresponding to the power injections, or more generally the \emph{drivers of uncertainty} in the system. A (finite dimensional) PCE basis corresponds to a set of $K$ polynomial basis functions $\Psib_k, \ k \in \mathcal{K} = \{0,\ldots, K-1\}$ such that
\begin{align}   \label{eq:orthogonality}
    \langle \Psib_l, \Psib_k \rangle = \mathbbm{E}\left[\Psib_l(\xib) \Psib_k(\xib) \right] = 0, \quad \mbox{for } l\neq k.
\end{align}
Each random variable $\mathbf{x} \in \{p_i,q_i,\vre_i,\vim_i\}$ in the system is expanded with respect to the basis functions
\begin{align}   \label{eq:pce_expansion}
    \mathbf{x} = \sum_{k \in \mathcal{K}} x_k \Psib_k(\xib),
\end{align}
where the scalars $x_k$ are the \emph{coefficients} of the PCE for $\mathbf{x}$. Uncertainty quantification then reduces to solving an extended system of power flow equations of the following form:
\begin{subequations} \label{eq:AC-PCE}
\begin{align}
    &\mbox{AC-PF equations},  \\
    &\mbox{PCE 1\textsuperscript{st} order PF equations} \\
    &\mbox{PCE 2\textsuperscript{nd} order PF equations} \\
    &\qquad \qquad \qquad   \vdots \nonumber
\end{align}
\end{subequations}
The details of these equations can be found in \cite{Muehlpfordt2019TPWRS} and are given in Table~\ref{tab:PCE_Reformulations} for completeness.

\begin{table*}
	\caption{Reformulations of power flow equations and moments in terms of \pce coefficients \cite{muhlpfordt2019chance}.}
	\vspace{-0.2in}
	\renewcommand{\arraystretch}{1.2}
	\label{tab:PCE_Reformulations}
	\begin{center}
		\begin{tabular}{l}
			\toprule
			Rectangular power flow in terms of \pce coefficients with $i \in \mathcal{N}$, $k \in \mathcal{K}$\\
			\hline
			$\langle \Psib_k {,} \Psib_k\rangle (p_{i,k})=  \sum_{j \in N} \sum_{k_1, k_2 \in \mathcal{K}} \langle \Psib_{k_1} \Psib_{k_2}, \Psib_{k} \rangle ( G_{ij} (\vre_{i,k_1} \vre_{j,k_2} + \vim_{i,k_1} \vim_{j,k_2}) + B_{ij} (\vim_{i,k_1} \vre_{j,k_2} - \vre_{i,k_1} \vim_{j,k_2}) ) $ \\
			$ \langle \Psib_k {,} \Psib_k\rangle (q_{i,k})=  \sum_{j \in N} \sum_{k_1, k_2 \in \mathcal{K}} \langle \Psib_{k_1} \Psib_{k_2}, \Psib_{k} \rangle ( G_{ij} ( \vim_{i,k_1} \vre_{j,k_2} - \vre_{i,k_1} \vim_{j,k_2}) - B_{ij} ( \vre_{i,k_1} \vre_{j,k_2} + \vim_{i,k_1} \vim_{j,k_2}) ) $ \\
			\hline
			Moments of squared line current magnitudes with $ij \in L$, $\vre_{ij,k} = \vre_{i,k} - \vre_{j,k}$, $\vim_{ij,k} = \vim_{i,k} - \vim_{j,k}$ \\
			$\mathbbm{E}\left[i_{i{\rightarrow}j}^2\right] = |y_{ij}^{\text{br}}|^2 \sum_{k \in \mathcal{K}}  \langle \Psib_k {,} \Psib_k\rangle ( (\vre_{ij,k} )^2 + ( \vim_{ij,k} )^2 )$  \\
			$ \sigma[ i_{i{\rightarrow}j}^2 ]^2 = |y_{ij}^{\text{br}}|^4 \sum_{k_1, k_2, k_3, k_4  \in \mathcal{K}} \langle \Psib_{k_1} \Psib_{k_2} \Psib_{k_3}{,} \Psib_{k_4}\rangle ( \vre_{i,k_1} \vre_{ij,k_2} \vre_{i,k_3} \vre_{ij,k_4} + 2 \vre_{ij,k_1} \vre_{ij,k_2} \vim_{ij,k_3} \vim_{ij,k_4} + \vim_{ij,k_1} \vim_{ij,k_2} \vim_{ij,k_3} \vim_{ij,k_4})  - \mathbbm{E}\left[i_{i\rightarrow j}^2\right]^2$ \\
			\hline
			Moments of squared voltage magnitudes with $i \in \mathcal{N}$ \\
			$\mathbbm{E}\left[V_{i}^2\right] = \sum_{k \in \mathcal{K}}  \langle \Psib_k {,} \Psib_k\rangle ( (\vre_{i,k})^2 + (\vim_{i,k})^2 )  $ \\
			$ \sigma[ V_{i}^2 ]^2 = \sum_{k_1, k_2, k_3, k_4  \in \mathcal{K}} \langle \Psib_{k_1} \Psib_{k_2} \Psib_{k_3}{,} \Psib_{k_4}\rangle ( \vre_{i,k_1} \vre_{i,k_2} \vre_{i,k_3} \vre_{i,k_4} + 2 \vre_{i,k_1} \vre_{i,k_2} \vim_{i,k_3} \vim_{i,k_4} + \vim_{i,k_1} \vim_{i,k_2} \vim_{i,k_3} \vim_{i,k_4})  - \mathbbm{E}\left[V_{i}^2\right]^2$ \\
			\bottomrule
		\end{tabular}
	\end{center}
\end{table*}

\subsection{Computational Complexity \& the Curse of Dimensionality} \label{subsec:curse}
A major advantage of PCE is that it reduces the infinite dimensional problem of uncertainty quantification into a finite dimensional problem. The accuracy of the resulting UQ depends on the number $K$ which denotes the number of basis functions used in the PCE expansion. This number depends on the \emph{degree} $\deg$ used in the expansion. Here $\deg$ denotes the maximum degree of the polynomial basis functions. As noted in \cite{Muehlpfordt2019TPWRS}, the number of basis functions grows exponentially in the chosen degree $\deg$ and is given by
\begin{align}   \label{eq:pce_scaling}
    |\mathcal{K}| = K = \dfrac{(n+\deg)!}{n!\deg!} \sim  \dfrac{n^{\deg}}{\deg!} \ \mbox{when } n\gg \deg. 
\end{align}
This phenomenon is a special case of the well-known \emph{curse of dimensionality}.
Fortunately, it was show in \cite{Muehlpfordt2019TPWRS} through several numerical studies that PCE with degree $2$ captures the non-linear nature of the PFE to a level of accuracy that is sufficient for all practical purposes. However, even for degree $2$, scaling the PCE method to large power system cases is computationally challenging. 

A closer inspection of the PCE-overloaded system of equations in Table~\ref{tab:PCE_Reformulations} and the scaling of $K$ in \eqref{eq:pce_scaling} shows that the computational complexity of the system of equations grows quite significantly with $K$. First, the number of equations and variables in the PCE overloaded system is $2NK$ compared to the $2N$ original PFEs in \eqref{eq:pfe}. Further, within each equation, expanding each of the quadratic terms and collecting the coefficients leads to a total of $O(K^2)$ terms. This level of scaling becomes quickly intractable for a large power system. As an example, consider a system with $N=1000$ buses and $n=10$ sources of uncertainty. To solve the PCE problem exactly with $\deg = 2$, we can compute $K = 66$ (using \eqref{eq:pce_scaling}) leading to $132,000$ equations with each equation having $K^2\sim 18,000$ terms! The $\deg = 1$ problem however, remains numerically tractable with $K_1 = 16$ and $K_1^2 = 256$.
However, as shown later $\deg = 2$ is needed to capture non trivial correlations effects between uncertainty sources that $\deg = 1$ does not capture. In this paper we develop approximation methods to reduce the computational complexity of $\deg=2$-PCE while still maintaining its high level of accuracy.

\section{\algo{}}   \label{sec:algorithm}

In this section we describe our approximations strategy that can significantly reduce the computational burden of solving the degree $2$ PCE-overloaded PFEs. Our approximations are based on two observations inherent to the system, (a) sparsity of the PCE coefficients, and (b) negligible contributions of higher order terms. These are described in detail in the subsections below.

\subsection{Sparsity of PCE Coefficients}    \label{subsec:sparsity}
The first key observation is that the PCE coefficients for all variables in the system are sparse, i.e, there are a large fraction of zero or near-zero coefficients. Among several explanations for why sparsity as a structural property might exist, a natural explanation -- re-inforced by our experimental observations, is that the value of a given quantity in the system (such as a given bus voltage) has strong dependence on only a few input power injections and essentially independent of the rest. Such structure can arise from factors such as geographical distance, where variables which are sufficiently far away from each other can be nearly independent. While there are several ways of discovering such independence properties, we use the PCE-overloaded PFE of degree $1$ for this purpose. This is descbribed in detail below.

In what follows, it will be useful to reformulate the PCE-overloaded PFE using a convenient matrix notation.
For a generic random variable $\mathbf{x}$, we can rewrite its PCE representation given in \eqref{eq:pce_expansion} as
\begin{align}  \label{eq:pce_matrix_form}
    \mathbf{x} = X^{(0)} + \langle X^{(1)},  \Psib^{(1)} \rangle + \langle X^{(2)},  \Psib^{(2)} \rangle,
\end{align}
where $X^{(0)}$ are the constant terms, the coefficient matrices $X^{(1)},X^{(2)}$ are defined as
\begin{align}  
    X^{(1)} = \begin{bmatrix}    
                    x_{1}\\
                    x_{2}\\
                    \vdots\\
                    x_{n}
                \end{bmatrix}, \quad 
    X^{(2)} = 
                \begin{bmatrix}    
                    x_{11} & x_{12} & \cdots & x_{1n}\\
                    x_{21} & x_{22} & \cdots & x_{1n}\\
                           &        & \vdots &      \\
                    x_{n1} & x_{n2} & \cdots & x_{nn}
                \end{bmatrix}.
\end{align}
The matrix of basis functions $\Psib^{(1)}$ contains all basis functions that are degree $1$ polynomials, and and $\Psib^{(2)}$ contains all degree $2$ basis functions re-arranged into a convenient matrix form such that $\Psib^{(2)}_{ij}(\xib)$ only depends on $\xi_i$ and $\xi_j$. 

We first solve the PCE-PFE with $\deg=1$ which amounts to setting $X^{(2)} = 0$, to obtain the coefficients $X^{(0)}_{\deg=1}$ and $X^{(1)}_{\deg=1}$. Whenever a coefficient $X_{j,\deg=1}^{(1)}$ is sufficiently small, we conclude that the random variable $\mathbf{x}$ is essentially independent of uncertainty source $\xi_i$. Further, since in the PCE expansion with $\deg=2$ the coefficient $X_{ij,\deg=2}$ corresponds to the basis function $\Psib^{(2)}_{ij}$ which is a function of only $\xi_i$ and $\xi_j$, we expect the near-independence of $\mathbf{x}$ on $\xi_i$ to be reflected in the degree $2$ coefficients via $X_{ij}^{(2)} \approx 0$. Based on this reasoning, we set a degree $2$ coefficient to $0$ \emph{apriori} if and only if
\begin{align} \label{eq:cutoff}
    |X_{i,\deg=1}^{(1)} X_{j,\deg=1}^{(1)}| < \Co \max_k |X_{k,\deg=1}^{(1)}| .
\end{align}

Forcing a fraction of the coefficients to zero will invariably lead to a reduction in accuracy. However, 
as we will show in Section~\ref{sec:numerical} through numerical experiments, using a well chosen cut-off $\Co$ one can reach a sparsity level of $15\% - 50\%$ with almost no loss of accuracy.

\subsection{Contribution of Higher Order Terms} \label{subsec:higher_order_terms}
Next, we seek to alleviate the large number of terms $O(K^2)$ in the quadratic expansion. Recall that this number was $\sim 18,000$ for $N=1000$ and $n=10$. This is achieved by observing that the contribution from large number of these terms that correspond to $4th$-order terms is negligible. These $4th$ order terms are generated by multiplying two degree $2$ basis functions in \eqref{eq:pce_expansion}. This reduces the number of terms drasitically from $O(K^2)$ to $O(K_1 K)$. 
In the previous example this reduces the number of terms from $\sim 18,000$ to $\sim 2,000$. As with the sparsification strategy in Section~\ref{subsec:sparsity}, we show through experiments that while using this approximation there is negligible loss of accuracy.

\subsection{Error Minimization and Warm Starting}    \label{subsec:overall_algorithm}
Finally, instead of solving a system of non-linear equations using e.g. Newton's method, we formulate an optimization problem that minimizes the error in the system of PCE-PFE. This step is necessary since, after setting a large fraction of coefficients to zero to enforce sparsity as described in Section~\ref{subsec:sparsity}, there are too few degrees of freedom, i.e., more constraints than variables. The cost function of the error-minimizing optimization problem can be chosen in several ways -- here we choose to use the $\ell_2$ loss function. Let $\mathbb{C}(\text{PCE coefficients}) = 0$ denote the set of constraints in Table~\ref{tab:PCE_Reformulations}. In the final step we solve the following unconstrained optimization problem:
\begin{align}   \label{eq:error_minimization}
    \min_{\text{PCE coefficients}} \|\mathbb{C}(\text{PCE coefficients})\|_2^2. 
\end{align}

Additionally, we choose a warm starting method to accelerate the optimization. We warm-start the PCE for $\deg=1$ with the solution to the deterministic PF, and the optimization problem in \eqref{eq:error_minimization}, with the coefficients obtained from the $\deg=1$ PCE. The steps of the algorithm are given in the pseudo-code.

\begin{algorithm}
\SetAlgoLined
Solve the $\deg=0$ PCE by solving the set of deterministic PFE ignoring the uncertainty
(For $N \sim 1,000$ is very fast $\sim\SI{0.1}{s}$)
Obtain $X^{(0)}_{\deg=0}$ - the deterministic power flow solution \;

 Solve the $\deg = 1$ AC-PF problem using the degree 0 result $x_0$ to warm start the $k=0$ coefficients. All degree $1$ coefficients are initialized to 0. Obtain $X^{(0)}_{\deg=1}$ and $X^{(1)}_{\deg=1}$, the order $0$ and order $1$ coefficients for PCE $\deg=1$. (For $N \sim 1,000$ is fast $\sim \SI{10}{s}$) \;
 
 Using $X^{(1)}_{\deg=1}$ and the thresholding procedure in \eqref{eq:cutoff}, determine a set $\mathcal{C}$ of degree $2$ coefficients that are expected to be small\;
 
 Set up the optimization problem in \eqref{eq:error_minimization} and \textbf{remove} the variables $X_i$ for $i \in \mathcal{C}$ from the problem. Solve \eqref{eq:error_minimization} to obtain the $\deg=2$ PCE solution $X^{(0)}_{\deg=2}, X^{(1)}_{\deg=2}, X^{(2)}_{\deg=2}$.

{\bf return} $X^{(0)}_{\deg=2}, X^{(1)}_{\deg=2}, X^{(2)}_{\deg=2}$\;
\caption{\algo{}}
\label{alg:overall_algorithm}
\end{algorithm}
\vspace{-0.1in}

\begin{figure}
    \includegraphics[width=0.4\textwidth]{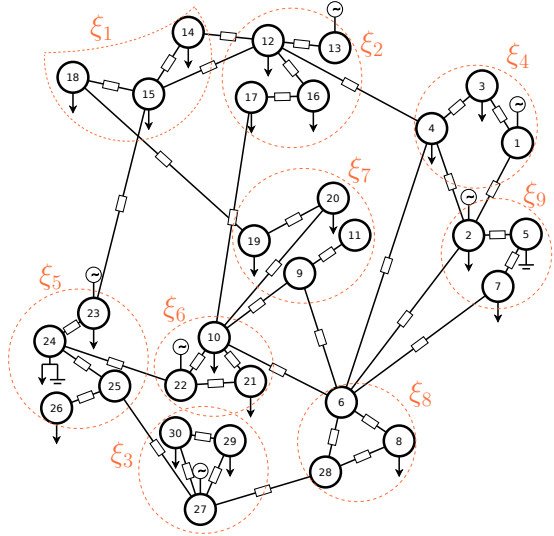}
    \caption{Illustration of the area decomposition of loads and fluctuations on a network with $30$ buses and $n = 9$ independent zones. Buses with down arrows are the (uncertain) loads.}
    \label{fig:case30}
    \vspace{-0.2in}
\end{figure}

\section{Numerical Simulations} \label{sec:numerical}

In this Section, we compare on a large network of 1354 buses the performances of \algo{}, our proposed UQ method, with respect to the standard PCE method and the Monte-Carlo method. We examine the computational time and the UQ accuracy of these techniques under two cases of load fluctuations: extreme and moderate fluctuations.

Computation are performed using computers with the same Intel Broadwell architecture and possessing 125GB of memory.
\subsection{Test-Cases Description}
The network that we consider is the test-case \verb| pglib_opf_case1354_pegase.m| contained in the IEEE PES Power Grid Library \cite{babaeinejadsarookolaee2019power}. This network features $1354$ buses and $673$ loads. The nominal values (without uncertainty) of our test-cases are set to be the default parameters of this network.

Uncertainty is produced by load fluctuations that depend on a particular geographical area. Loads in the network are partitioned in $n$ geographical areas based on their network proximity using the recursive graph partitioning method of the METIS package \cite{metis}. Each area is associated with an independent centered and normalized random variables $\xi_a$ that represents the type of uncertainty in the area $a$, see Figure~\ref{fig:case30}. The load fluctuations within each zone are fully correlated, and between two different zones are independent, i.e., for two different zones $a$ and $a'$ $\E(\xi_a \xi_{a'}) =0$.
Lastly the intensity of fluctuations of loads is controlled by a parameter $\epsilon$ such that the active and reactive power of a single load in an area $a$ reads as follows,
\begin{align}
            p_{\text{load}} = p^{\text{nominal}}_{\text{load}} (1 + \xi_a \epsilon), \quad  q_{\text{load}} = q^{\text{nominal}}_{\text{load}} (1 + \xi_a \epsilon ).\label{eq:power_at_load}
\end{align}
Note that Eq.~\eqref{eq:power_at_load} implies that the nominal power factor of each load is kept constant while its total power consumption varies by an amount proportional to $\epsilon$.
Generators are assumed to respond to active power fluctuations uniformly i.e. 
\begin{align}   \label{eq:recourse}
            p_{\text{gen}} = p_{\text{gen}}^{\text{nom}} + \frac{P_{\text{load}}-P^{\text{nom}}_{\text{load}}}{N_{\text{gen}}}, \quad V_{\text{gen}} = V_{\text{gen}}^{\text{nom}}, 
\end{align}
where $P^{\text{nom}}_{\text{load}}$ and $P_{\text{load}}$ are the total nominal and realized load consumption respectively and $N_{\text{gen}}$ is the total number of generators. The recourse policy in \eqref{eq:recourse} corresponds to automatic generation control (AGC) with uniform participation factors. Other policies can be incorporated in a similar way. 

The two different cases of load fluctuations that we consider are obtain through varying $\epsilon$.
The {\bf extreme} case corresponds to power fluctuations of $3\%$ at each of the 673 loads in the network. This is implemented by setting $\epsilon = 3\times 10^{-2}$. The {\bf moderate} case corresponds to power fluctuations of $1\%$ at each load in the entire network by setting $\epsilon = 10^{-2}$. For these two cases, the random variable $\xi_a$ modelling the uncertainty per area are chosen to be normalized uniform distributions.

\subsection{Choice of Hyperparameter $\Co$}
As mentioned in Subsection~\ref{subsec:sparsity}, \algo{} depends on a tunable hyperparameter $\Co$ that promotes sparsity of degree 2 coefficients based on the computation of degree 1 coefficients, see Eq.~\eqref{eq:cutoff}. 
Computational time and degree 2 coefficient sparsity for two different values of the cutoff $\Co=0$ and $\Co=10^{-10}$ are displayed in Table~\ref{tab:cutoff_results}. These results are obtained on the 1354 buses test-case with 10 areas and extreme fluctuations.
\begin{table}[]
    \centering
    \begin{tabular}{l|c c c c c}
        \\
        $\Co$ & $0$ & $10^{-10}$  \\
        \hline\\
        Computational Time (\si{s})  &  647 &  522  \\
        Degree 2 Coefficients Sparsity  \% & 0& 39
    \end{tabular}
    \caption{Computational time and sparsity of \algo{} for two different value of the cutoff hyperparameter.}
    \label{tab:cutoff_results}
\vspace{-0.3in}
\end{table}

We see that a very small cutoff of $10^{-10}$ is already sufficient to remove around $40\%$ of the PCE degree 2 coefficients which leads to an overall speed-up of $20\%$. Moreover with a cutoff of only $10^{-10}$, the impacts on the UQ quality remains unnoticeable.
For all numerical simulations, we chose the cutoff of \algo{} to be $\Co=10^{-10}$ which in practice leads to a good trade-off between sparsity and accuracy.

\subsection{Accuracy of Uncertainty Quantification}\label{subsec:accuracy_uq}
\begin{figure}
\includegraphics[width=0.45\textwidth]{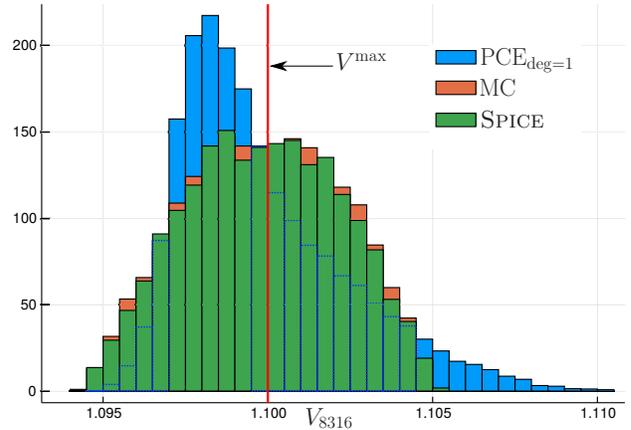}
    \caption{Voltage distribution at bus label 8316 computed on the 1354 buses test-case with 10 areas and extreme fluctuations using $10^4$ samples. Histograms obtained using Monte-Carlo is in orange, full PCE of degree 1 in blue and \algo{} in green. The voltage limit for this bus is displayed in red.
    }
    \label{fig:histo}
    \vspace{-0.1in}
\end{figure}

The output of an UQ method, whether it runs with Monte-Carlo or PCE, is a probability distribution for each variable in the system given in the form of a histogram. These histograms are produced by first drawing $M=10^4$ realizations of the random variables $\xi_a$ and then for each of them, either solve a Power-Flow problem if one uses Monte-Carlo or evaluate a polynomial if one uses a PCE based approach. Finally, the output of the PF problems or polynomial evaluations are aggregated into discrete histograms using bins of size $5 \times 10^{-3}$ times the typical variable scale computed from the system bounds and chosen to be $V_{\text{max}} - V_{\text{min}}$ for voltages and $S^{\text{max}}_{i \rightarrow j}$ for line power flows. Note that a bin size much smaller than $\approx M^{-1/2}$ goes beyond the precision that one would expect to achieve using $M$ samples. Typical histograms obtained with this procedure are displayed in Figure~\ref{fig:histo}.

The distance between two histograms $h_1$ and $h_2$ is measured using the Total Variation (TV) distance,
\begin{equation}
    \TV(h_1,h_2) = \frac{1}{2M}\sum_{b\in\text{Bins}} |h_1(b) - h_2(b)|, \label{eq:TV_def}
\end{equation}
which accounts for the average difference in counts for each histogram. The reason for choosing this metric is that is translates directly into guarantees for computing probabilities: The difference in the probability of events computed from two histograms with TV distance of $\delta$ is no larger than $\delta$.

Average and maximum TV distance between histograms of voltage magnitude and power line flows are presented in Table~\ref{tab:tv_results}. The results are obtained on the 1354 test-case with 10 areas and extreme fluctuations. The four methods are the full PCE of degree 1 and degree 2, \algo{} and Monte-Carlo. The TV distance of each method is measured with respect to a reference histogram produced by Monte-Carlo using an independent draw of $10^4$ samples.
\begin{table}[]
    \centering
    \begin{tabular}{l|c c c c c}
        \\
        Method & $\deg = 1$ & $\deg = 2$& \algo{} & MC \\
        \hline\\
        Computational Time (\si{s})  &  6 & 1714  & 522 & 1060\\
        Ave. TV, Voltage Magnitude & 0.10 & 0.006  & 0.005 &  0.005   \\
        Max. TV, Voltage Magnitude & 0.80 & 0.031  & 0.029 &  0.033 \\
        Ave. TV, Power Line Flow & 0.015 & 0.012  & 0.012 &  0.012   \\
        Max. TV, Power Line Flow & 0.19 & 0.047  & 0.043 &  0.050 \\
    \end{tabular}
    \caption{Average and maximum TV distance between UQ methods with respect to a reference Monte-Carlo histogram.}
    \label{tab:tv_results}
    \vspace{-0.3in}
\end{table}

Note that the TV distance between two independent Monte-Carlo runs is not zero but is about $0.5\%$ to $1.0\%$ on average and between $3\%$ to $5\%$ in the worse case. These discrepancies are caused by unavoidable statistical fluctuations that arose in our finite sample set. Therefore, the results for Monte-Carlo should be seen as a reference for the minimum TV distance we can expect to achieve with $M=10^{4}$ samples.
The quality of UQ for PCE of degree 2 and \algo{} are similar and indistinguishable from Monte-Carlo, while \algo{} being much less computationally intensive. PCE of degree 1 had the advantage to be extremely fast to run, but it performs poorly on this test-case with extreme load fluctuations of $3\%$. We obtained for voltage magnitude a TV distance of $10\%$ on average and it goes up to $80\%$ in the worst case.
These results show that PCE of degree 2 is a well-suited method for uncertainty quantification for PF equations. Moreover, the sparsity promoting techniques implemented in $\algo{}$ does not come at a noticeable cost in terms of accuracy.

\subsection{Uncertainty Quantification Robustness with \algo{}}
We test how robust is the UQ accuracy of the PCE method to a change in the distribution of load fluctuations. The PCE coefficients are computed with \algo{} using the 1354 buses test-case with extreme fluctuations. This setting is identical to what is described in the previous Subsection~\ref{subsec:accuracy_uq}, for which the random variables $\xi_a$ are normalized and centered \emph{uniform} distributions. However, unlike in Subsection~\ref{subsec:accuracy_uq}, the histograms are produced from the PCE polynomial evaluations using $M=10^4$ samples generated by variables $\xi_a$ that are chosen to be normalized and centered \emph{Gaussian} distributions. We compare \algo{} with two Monte-Carlo runs, one for which the $M=10^4$ samples arose from a Gaussian distributions and one for which the samples come from the uniform distribution. 

Results of average  and  maximum  TV  distance  between  histograms of  voltage  magnitude  and  power  line  flows  are  presented in Table~\ref{tab:robust_tv} for \algo{} and Monte-Carlo. The TV distance is measured with respect to a reference Monte-Carlo histogram produced from an independent draws of $10^4$ samples from Gaussian distributions. We see that the Gaussian fluctuations produce very different histograms than in the uniform case. The TV difference between the uniform and Gaussian Monte-Carlo for voltage magnitude is about $50\%$ on average and $100\%$ in the worse case, which means that there is no overlap at all between the two histograms.
The results also show that even though the PCE coefficients found by \algo{} are suited for a uniform distribution, it remains accurate when the uncertainty arises from a very different distribution. This highlights an important feature of PCE that, unlike Monte-Carlo, it computes a deterministic mapping between load fluctuations and the power flow variables of the system. This mapping can latter be reused with a different uncertainty source such as historical data or uncertainty scenarios at no extra cost and with little impact on the UQ quality.  

\begin{table}[]
    \centering
    \begin{tabular}{l|c c c}
        \\
        Method & \algo{} & $\text{MC}_{\text{Gaussian}}$ & $\text{MC}_{\text{uniform}}$\\
        \hline\\
        Ave. TV, Voltage Magnitude &  0.008 & 0.008 &  0.53 \\
        Max. TV, Voltage Magnitude &  0.037 & 0.040 & 1.00\\
        Ave. TV, Power Line Flow  & 0.015 &  0.015 &  0.12 \\
        Max. TV, Power Line Flow  & 0.054 &  0.050 & 1.00\\
    \end{tabular}
    \caption{TV distances for Gaussian fluctuations when PCE coefficients are  computed for uniform fluctuations.  
    } 
    \label{tab:robust_tv}
    \vspace{-0.3in}
\end{table}

\subsection{Computational Speed}
\begin{figure}[h]
    \includegraphics[width=0.5\textwidth]{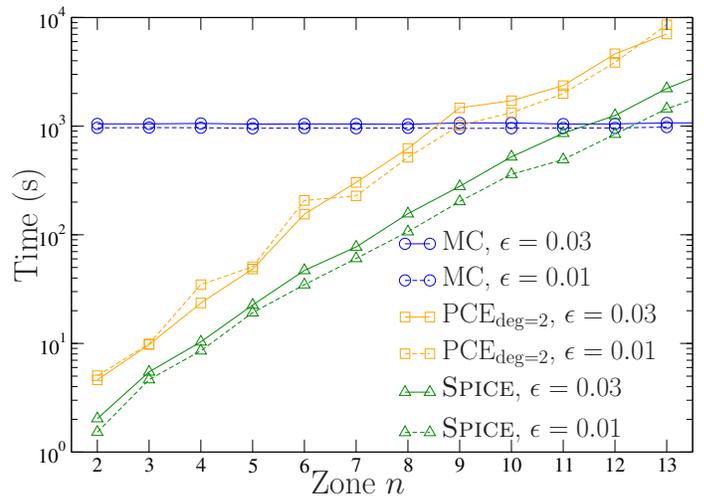}
    \caption{Computational time taken by three different UQ methods on the 1354 buses test-case for different number $n$ zones. Computational times are shown in blue, yellow and green for Monte-Carlo, Full PCE of degree 2 and \algo{} respectively. Solid lines are associated with extreme load fluctuations of $3\%$ and dashed lines are associated with moderate load fluctuations of $1\%$. On average \algo{} is $3.5$ faster than the full PCE of degree 2}
    \label{fig:ipopt_times}
    \vspace{-0.2in}
\end{figure}

The computational speed of Monte-Carlo, Full PCE of degree 2 and $\algo{}$ are compared on the 1354 buses test-case for different number of areas ranging from $n=2$ to $n=13$ under both extreme and moderate fluctuations. The computational times of the different method are displayed in Figure~\ref{fig:ipopt_times}. We would like to stress that the times reported on that figure consists solely on computations necessary to produce the histograms (solving PF equations and PCE equations). In particular it does not account for the overhead time spent on saving and handling the datasets produced for each method. We will see in Subsection~\ref{subsec:memory} that this overhead is not negligible for Monte-Carlo methods and ends up multiplying the whole run-time by a factor 2 to 4.

As expected the Monte-Carlo methods are not sensitive to the number of sources of uncertainty present in the system and depends only on the time required to solve $10^4$ PF equations. For PCE techniques, the curse of dimensionality is apparent as the computational time increases exponentially with the number of zones. The benefits of the computational enhancements of \algo{} translates in average into a $3.5$ time speed-up with respect to the standard full PCE of degree 2. This makes \algo{} more competitive than Monte-Carlo even for large systems when the number of uncertainty sources is around 10. Note also that the computational time for the standard PCE is similar for extreme and moderate load fluctuactions while there is a $40\%$ time difference for \algo{}. This reason it that \algo{} takes advantage of the sparsity in the PCE coefficients. When the fluctuactions are moderate, loads have a lesser impact on variables located further away which leads to sparser PCE coefficients. 

\subsection{Overhead Time and Memory Storage}\label{subsec:memory}

As mentioned previously, the computational time and memory capacity required to handle and store the datasets produced by the UQ methods differs significantly between Monte-Carlo and \algo{}.

Concerning Monte-Carlo, one can only store the final histograms that are composed of $M=10^4$ points for every variable (voltage, active power, reactive power, line flows) that are at each of the 1354 buses. This ends up constituting a file of 0.6GB and multiplies the whole run-time of the algorithm by a factor 2 to 4 owing to data loading latencies. While typical computational times reported in Figure~\ref{fig:ipopt_times} are around 1000 seconds for Monte-Carlo, the whole run-time including storage and data handling reaches in practice 1 hour.

The story is different for PCE methods like \algo{} as it offers the capability to store only the non-zero PCE coefficients and generate the histograms later on the fly. The number of non-zero coefficients required for \algo{} is not more than a hundred per variable and per bus and can be stored using only a dozen of MB. Histograms are generated from evaluations of PCE polynomial with independent draws of the random variables $\xi_a$. Moreover, the evaluation of the PCE polynomials can be done efficiently using sparse matrix multiplication. For the 1354 buses system using 13 areas, the whole operation only takes 5 seconds. 

\section{Application to AC-OPF with Chance Constrains}\label{sec:CC-OPF}
In this section we apply our proposed UQ method \algo{} for solving stochastic AC-OPF with chance-constraints (CC-AC-OPF).
In this setting the power flow equations described in Section~\ref{sec:PFE} are supplemented with the traditional constraints on voltage, line power and generation limits enforced probabilistically (the so-called chance constraints),
\begin{align*}
    &\mathbb{P}\left(V_{i}^{\rm{min}}\leq V_{i} \leq V_{i}^{\rm{max}}\right)\geq 1 - \delta,  &\text{Voltage Limit}\\
    &\mathbb{P}\left(S_{i\to j}^2\leq (S_{i\to j}^{\rm{max}})^2\right)\geq 1 - \delta,  &\text{Power Line Limit}\\
    &\mathbb{P}\left(p_i^{\rm{min}} \leq p_i\leq p_i^{\rm{max}}\right)\geq 1 - \delta,  &\text{Active Generation Limit}\\
    &\mathbb{P}\left(q_i^{\rm{min}} \leq q_i\leq q_i^{\rm{max}}\right)\geq 1 - \delta,  &\text{Reactive Generation Limit}
\end{align*}
where the confidence level at which each constraint are satisfied is $1 - \delta$.  
The procedure that we implement for solving the CC-AC-OPF problem is described by Algorithm~\ref{alg:CC_OPF_algorithm}. It is an iterative scheme that goes back and forth between solving a \emph{deterministic} AC-OPF problem with \emph{effective} voltage, power line and generation limits and a UQ evaluation of the chance constraints with \algo{} to update the effective bounds, see \cite{Roald18} for more details.

\begin{algorithm}
\SetAlgoLined
Initialization of effective voltage, power and generation limits:
$V^{\rm{min}}_{\rm eff},V^{\rm{max}}_{\rm eff},S^{\rm{max}}_{\rm eff},\ldots \longleftarrow V^{\rm{min}},V^{\rm{max}},S^{\rm{max}},\ldots$ \;

\Repeat{all excess differences vanishes $\Delta V^{\rm{min}}=\Delta V^{\rm{max}}=\Delta S^{\rm{max}}=\ldots=0$}{
Run deterministic AC-OPF with effective limits to determine operating point\;
At the current operating point, evaluate with \algo{} the $\delta$ quantiles $Q_{V^{\rm{min}}}, Q_{V^{\rm{max}}}, Q_{S^{\rm{max}}}, \ldots$ of each chance constraints (e.g. $\mathbb{P}\left(V\leq Q_{V^{\rm{max}}}\right) = 1-\delta $)\;
Compute excess differences between limits and quantiles:
$\Delta V^{\rm{min}} \longleftarrow \min(Q_{V^{\rm{min}}} - V^{\rm{min}}, 0)$,  
$\Delta V^{\rm{max}} \longleftarrow \max(Q_{V^{\rm{max}}} - V^{\rm{max}}, 0)$,
\ldots\;

Update effective limits with excess differences:
$V^{\rm{min}}_{\rm eff} \longleftarrow V^{\rm{min}}_{\rm eff} + \Delta V^{\rm{min}}$,  
$V^{\rm{max}}_{\rm eff} \longleftarrow V^{\rm{max}}_{\rm eff} + \Delta V^{\rm{max}}$,
\ldots\;
}

\caption{Iterative CC-AC-OPF with \algo{}}
\label{alg:CC_OPF_algorithm}
\end{algorithm}
\vspace{-0.1in}

Our test-case for CC-AC-OPF is the 1354 buses system described in Section~\ref{sec:numerical} with moderate fluctuations and $n=10$ areas.
Bounds on the reactive power at generator 46 located at bus label 1754 have been removed as they were too restrictive for admitting a feasible solution to the stochastic CC-AC-OPF with $1\%$ load fluctuations.

\begin{table}[]
    \centering
    \begin{tabular}{l|c c c}
        
        CC-AC-OPF confidence level   $(1-\delta)$   &  $95\%$ &  $99\%$\\
        \hline\\
        Number of Iterations & 4 & 5 \\
        Computational Time (\si{s})  &  1463 & 1793\\
        Post-Validation with MC & \checkmark & \checkmark
    \end{tabular}
    \caption{Number of iterations and computational time for solving CC-AC-OPF with Algorithm~\ref{alg:CC_OPF_algorithm}.} 
    \label{tab:cc-ac-opf}
    \vspace{-0.3in}
\end{table}

Results on convergence time and number of iterations necessary to solve CC-AC-OPF with \algo{} are shown in Table~\ref{tab:cc-ac-opf}. Once the optimal solution is returned by Algorithm~\ref{alg:CC_OPF_algorithm}, the probability of bound violation are verified using 3 independent Monte-Carlo validations.
We have also tested Algorithm~\ref{alg:CC_OPF_algorithm} using PCE of degree 1 instead of \algo{}. However for a confidence level of $99\%$, the solution provided by PCE of degree 1 underestimates the reactive power fluctuations arising at generator 16 (label 757) which ends up violating its limit for more than $1\%$.

\section{Conclusion and Future Work}
In this paper, we have proposed an efficient and accurate UQ method, \algo{}, for characterizing uncertainty in AC power flow equations which is a computationally enhanced PCE method of order 2.
The main advantages of \algo{} are that a) it scales to large systems and is computationally superior compared to Monte-Carlo b) it takes advantage of the inherent sparsity pattern of the fluctuation responses c) it is robust with respect to changes in the uncertainty distribution and d) it requires a low amount of memory for data storage.

In the future, we will focus our effort on developing a tractable single optimization formulation that incorporates \algo{} within the AC-OPF problem directly. This will overcome the need for going through the iterative Algorithm~\ref{alg:CC_OPF_algorithm} for solving the CC-AC-OPF problem and potentially increase the speed by another factor 4-5.

\bibliographystyle{IEEEtran}
\bibliography{IEEEabrv}

\end{document}